\documentclass[12pt]{article}
\usepackage{graphicx}
\usepackage{indentfirst}

\def\P{{\bf P}}

\def\checkpage#1{
\dimen5 = \pagetotal
\ifdim \dimen5 < \pagegoal
\advance \dimen5 by #1
\ifdim \dimen5 > \pagegoal
\eject
\fi
\fi
}

\begin{document}

\baselineskip 18pt

\begin{center}
{\Large Optimal Strategies and Rules for the Game of Horse} \\
\bigskip
by \\
\medskip
Daniel Rosenthal \quad and \quad Je{f}frey S.\ Rosenthal \\
\medskip
University of Toronto, Canada \\
\medskip
(December 2021; last revised January 9, 2022.) \\
\medskip
{\it Notices of the American Mathematical Society}, to appear.
\end{center}

\section*{Introduction}

The game ``Horse'' is a basketball variant in which two players attempt
various basketball shots of their choosing, each trying to make shots
which the other cannot.  It is a popular past-time on basketball courts
around the world, and has even been played by professional players
(NBA.com, 2020).

Over the years, we have noticed that some players
(including our own father or ``Pops'')
tended to choose very {\it
easy}\/ shots, which seemed to work to their advantage.
This made us wonder
if we could analyse, mathematically, the effect of shot difficulty on the
probability of winning points when playing Horse.  In this paper, we will
use simple probability calculations to show that, under the Traditional
Rules (TR), it is indeed often optimal to choose very easy shots.  We will
also introduce modified rules which we call the Pops Rules
(PR), and show that these modified rules lead to the different (and, we
think, better) situation in which it is optimal to choose more
difficult shots.

We shall first assume that two players of equal ability are playing Horse.
Player One can select a shot which each player has the same
probability $p$ of making, for any $p\in(0,1)$.  We will consider the
following two questions:

Q1: As a function of $p$, what is the probability (under each of TR and PR)
that Player One will score a point on their turn?

Q2: What choice of $p$ maximises the probability (under each of TR and PR)
that Player One will score a point on their turn?

We will then discuss what these calculations tell us about which rules,
TR or PR, are preferable when playing Horse.

Finally, we will consider more general situations in which the two
players have unequal abilities.

There have been a few previous statistical analyses of Horse.
Tarpey and Ogden (2016) use logistic regression to model their own
shot success probabilities as a function of
distance to the basketball rim, and then model the game of Horse
as a Markov chain in an effort to optimise win probabilities.
And, Stevenson (2020) presents a statistical model for shot selection in
a simulated hypothetical Horse match between two specific NBA players.

\section*{The Game of Horse}

The game Horse is played in turns, as follows.
Let's suppose it is Player One's turn.
Player One first describes and attempts a basketball shot of their choosing.
If they miss the shot, then their turn is over and it is then Player Two's
turn.
But if Player One makes their shot,
then Player Two has to attempt the same shot,
and if Player Two misses that shot
then they receive a penalty of the next letter
in the word H-O-R-S-E (i.e.\ Player One is awarded a point).
If Player One makes their shot and Player Two also makes
the shot, then no points are awarded.
The game ends when one player has received all of the letters H-O-R-S-E,
at which point that player loses.

In summary, there are three possible scenarios for
Player One's turn in a game of Horse:

S1. Player One fails to make their shot.
In this case, no points are awarded, and Player Two's turn begins.

S2. Player One succeeds in making their shot,
and then Player Two fails to make that same shot.
In this case, Player Two receives a letter (i.e.\ Player One scores a point).
(Player One then restarts their turn,
and might score additional points, though that does not 
not affect our considerations herein.)

S3. Player One succeeds in making their shot,
and then Player Two succeeds in making that same shot.

Under scenario S3, again no points are awarded.  But what happens next?
In the Traditional Rules (TR), after S3, Player One's turn restarts,
i.e.\ Player One then describes and attempts a new shot of their
choosing and the game proceeds according to the above three scenarios.  We
shall also consider alternative rules, which we call the Pop Rules (PR),
wherein scenario S3 ends Player One's turn, i.e.\ the turn
automatically ``pops'' over to Player Two after both players
successfully make the same shot.

\section*{Q1: The Probability of Scoring a Point}

We assume now that Player One selects a shot which each player has
the same probability~$p \in (0,1)$ of making.
Then the probabilities for
the above three scenarios are easily seen to be:

\def\fail{{\rm fail}}
\def\succeed{{\rm succeed}}

$\P(S1) \ = \ \P\big(\hbox{Player One fails}\big)
\ = \ 1-p$.

$\P(S2) \ = \ \P\big(\hbox{Player One succeeds, then Player Two fails}\big)
\ = \ p \, (1-p)$.

$\P(S3) \ = \ \P\big(\hbox{Player One succeeds, then Player Two succeeds}\big)
\ = \ p^2$.

\medskip

Let $A(p)$ be the probability that Player One will score
a point on their first turn under TR,
and let $B(p)$ be the same probability under PR.

Then clearly $B(p) = \P(S2) = p(1-p)$.

On the other hand, since under TR after scenario S3 the turn restarts,
$$
A(p) \ = \ \P(S2) + \P(S3) \, A(p) \ = \ p(1-p) + p^2 \, A(p)
\, .
$$
So, solving this equation for $A(p)$ gives
$$
A(p)
\ = \ {p(1-p) \over 1-p^2}
\ = \ {p(1-p) \over (1-p)(1+p)}
\ = \ {p \over 1+p}
\, .
$$

\begin{figure}[ht]
\centerline{\includegraphics[width=10cm]{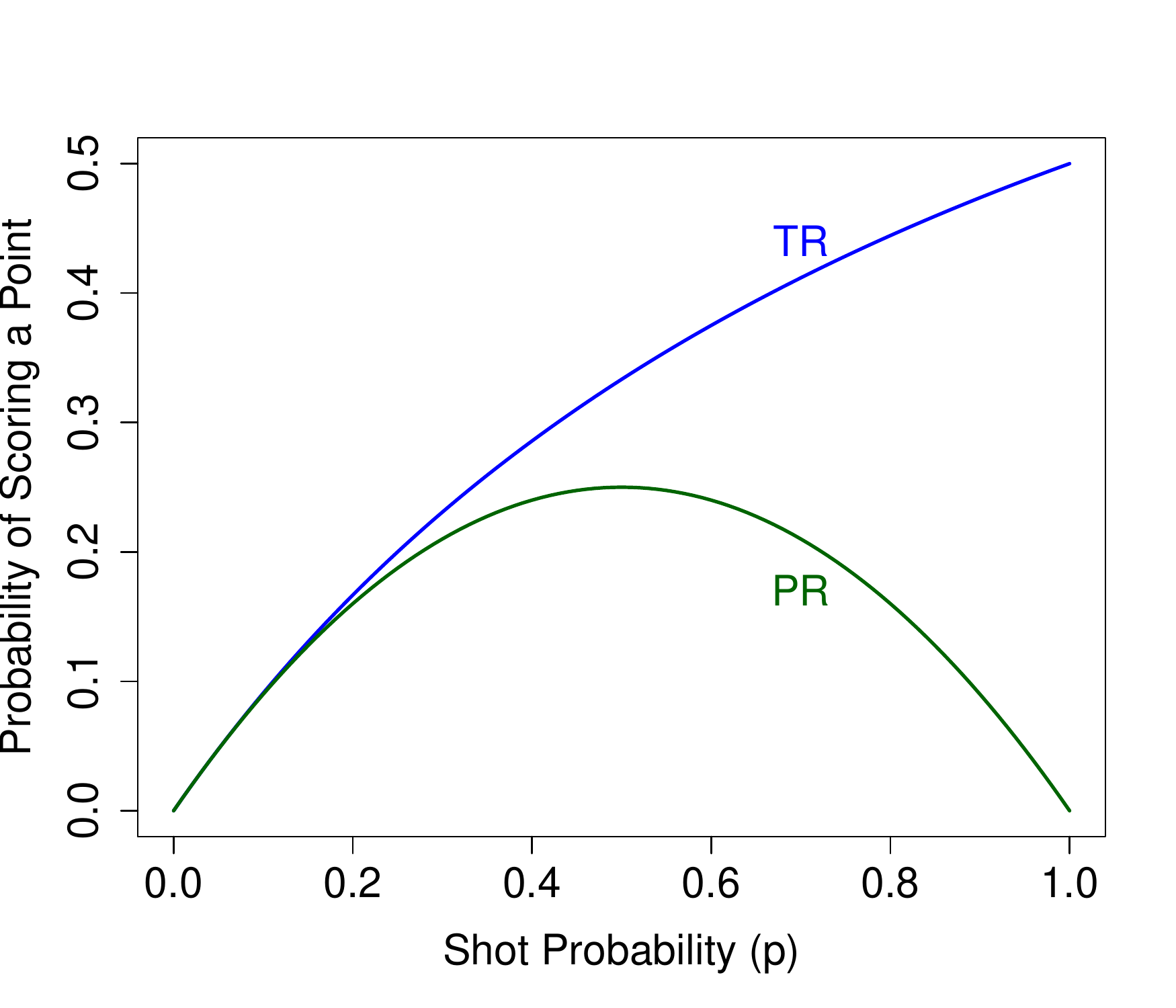}}
\caption{\sl
Graph of the probability of Player One
scoring a point in Horse, as a function of the
shot probability $p$, under the Traditional Rules (TR: blue, top)
and under the Pops Rules (PR: green, bottom).
}
\label{fig-equal}
\end{figure}

Graphs of the functions $A(p)$ and $B(p)$ are shown in
Figure~\ref{fig-equal}.  Note in particular that we always have $A(p) \ge
B(p)$, i.e.\ the restarting of the turn after S3 under the TR
can only {\it help} Player One to score.

\section*{Q2: The Optimal Success Probability $p$}

Under PR, since the probability of scoring a point is
$B(p) = p(1-p)$, it
is maximised when $p=1/2$.  That is, under PR it is optimal for the first
player to choose a medium-level shot with probability 50\% of success.

However, under TR, the probability of scoring a point is
$A(p) = {p \over 1+p}$ which is an
increasing function of $p$, so it is maximised as $p \nearrow 1$.  That
is, under TR it is optimal for Player One to choose an extremely
easy shot with probability of success near 100\%.  (We exclude the
case where $p$ is actually equal to~1, as being inconceivable and also
leading to a turn that never ends.)

Intuitively, this optimal TR play will lead to a very long turn in which
each player has approximately equal probability of being the first to
finally fail, leading to a probability of nearly 50\% of the first
player scoring a point.

\section*{Discussion: The Optimal Rule Choice}

The above results show that, if Horse is played optimally under the
Traditional Rules (TR), then players are motivated to take extremely easy
shots.  This will lead to a very uninteresting game, in which only easy
and boring shots are attempted, and nearly all shots are made, and turns
last an extremely long time until a point is finally scored.

By contrast, if Horse is played optimally under the modified Pops Rules (PR),
then players are motivated to take shots which have success probability
near 50\%.  This will lead to an interesting game, in which players take
interesting shots, and succeed and fail approximately equally often.

We conclude from our analysis that when playing Horse,
the modified Pops Rules (PR) lead to more interesting optimal strategies
and thus a better game.  We therefore believe that the
Pops Rules should be used instead of the Traditional Rules.

\section*{Unequal Players}

So far we have assumed that each player has the same probability of making
each shot.  But now we will consider the case where the probabilities are
not equal.  Suppose the players have probabilities $p_1$ and $p_2$,
respectively, of making a shot.
Then the above scenario probabilities become

$\P(S1) \ = \ \P\big(\hbox{Player One fails}\big)
\ = \ 1-p_1$.

$\P(S2) \ = \ \P\big(\hbox{Player One succeeds, then Player Two fails}\big)
\ = \ p_1 \, (1-p_2)$.

$\P(S3) \ = \ \P\big(\hbox{Player One succeeds, then Player Two succeeds}\big)
\ = \ p_1 \, p_2$.

\medskip

\begin{figure}[ht]
\centerline{\includegraphics[width=14cm]{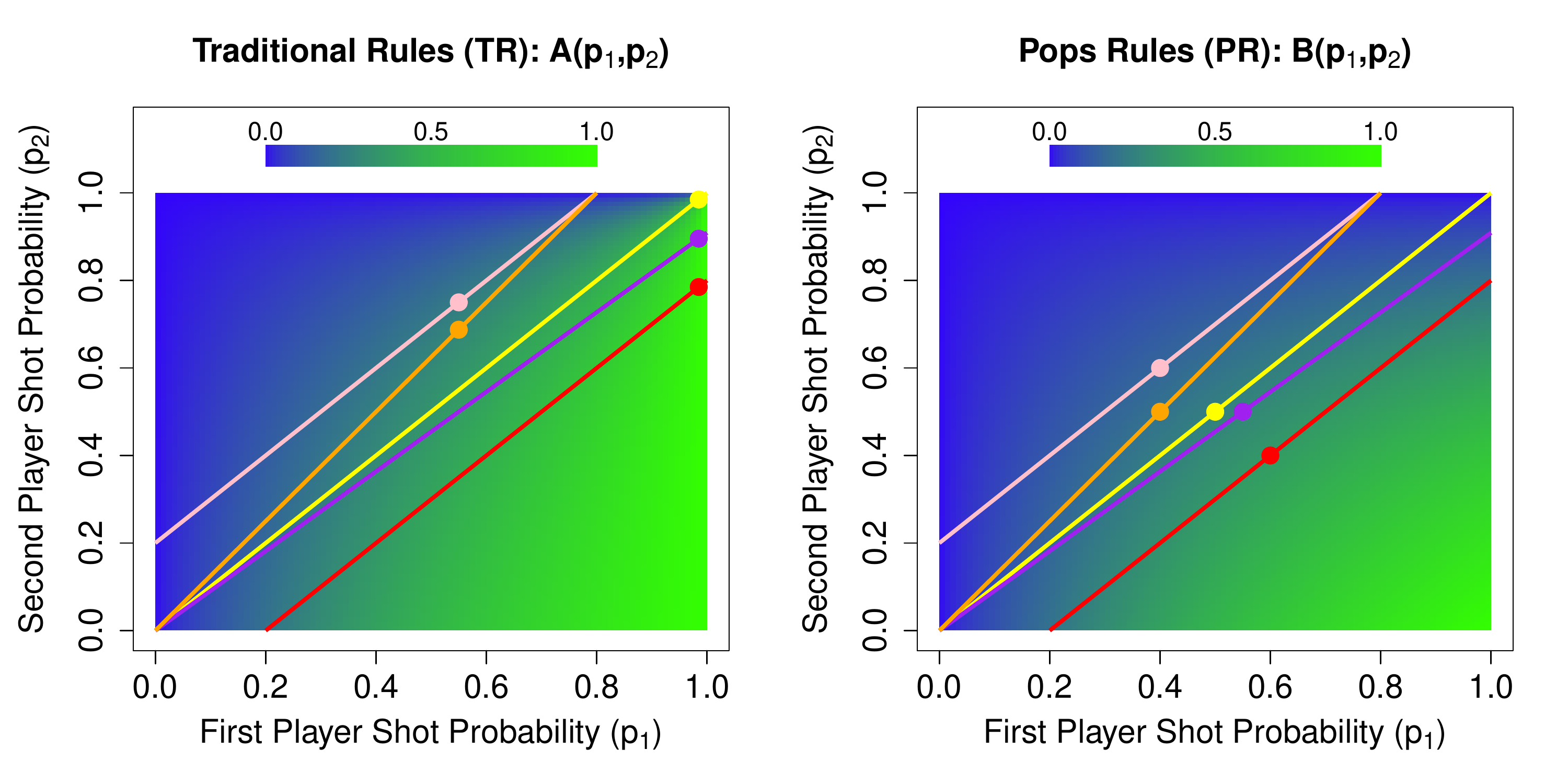}}
\caption{\sl
Heat map of the probability of Player One
scoring a point on their turn in a game of Horse, as a function of the two
players' shot probabilities $p_1$ and $p_2$, under the Traditional Rules (left)
and under the Pops Rules (right).
Also shown (from top to bottom) are the lines where
$p_2=p_1+0.2$ (pink),
$p_2=p_1/0.8$ (orange),
$p_2=p_1$ (yellow),
$p_2=p_1/1.1$ (purple),
and $p_2=p_1-0.2$ (red),
together with the optimal $p_1$ in each case.
}
\label{fig-unequal}
\end{figure}

Let $A(p_1,p_2)$ and $B(p_1,p_2)$ be
the probability that Player One will score
a point on their first turn given $p_1$ and $p_2$
under TR and under PR, respectively.
Then, similar to the above, under PR we have $B(p_1,p_2) = \P(S2) = p_1(1-p_2)$,
and under TR we have
$$
A(p_1,p_2) \ = \ \P(S2) + \P(S3) \, A(p_1,p_2) \ = \ p_1 \, (1-p_2) + p_1 \, p_2 \, A(p_1,p_2)
\, .
$$
So, solving this equation for $A(p_1,p_2)$ gives
$$
A(p_1,p_2)
\ = \ {p_1 \, (1-p_2) \over 1 - p_1 \, p_2}
\, .
$$
A heatmap of
the functions $A(p_1,p_2)$ and $B(p_1,p_2)$ are shown in Figure~\ref{fig-unequal}.
Of course, these probabilities are highest when $p_1$ is large and $p_2$ is small.
Once again, we always have $A(p_1,p_2) \ge B(p_1,p_2)$ as we must.

\checkpage{2cm}
\bigskip\noindent\bf Linear Difference case: \rm

Consider the linear difference case, where $p_2=p_1-d$ for
some fixed constant $d\in(-1,1)$.
To ensure that $p_1,p_2\in(0,1)$, we restrict to $p_1 \in (0, \min[1,1+d])$.
Then under PR,
\[
B(p_1,p_1-d)
\ = \ p_1(1-[p_1-d])
\ = \ (1+d) p_1 - p_1^2
\, .
\]
It follows that $B(p_1,p_1-d)$ is maximised when $p_1 = (1+d)/2$.  This holds
for any $d \in (-1,1)$, i.e.\ it is true regardless of which player is
better.

However, under TR,
\[
A(p_1,p_1-d)
\ = \ \frac{p_1(1-[p_1-d])}{1 - p_1[p_1-d]}
\ = \ \frac{p_1 - p_1^2 + p_1d}{1 - p_1^2 + p_1d}
\, .
\]
To analyse this, we compute that for fixed $d$,
$$
{d \over dp_1} \, A(p_1,p_1-d)
\ = \ {(1-p_1)^2 + d \over (1+dp_1-p_1^2)^2}
\, .
$$
If $d \ge 0$ (i.e.\ Player One is better than Player Two),
then ${d \over dp_1} \, A(p_1,p_1-d) > 0$ for all $p_1<1$, i.e.\
$A(p_1,p_1-d)$ is again an increasing function of $p_1$,
so it is again maximised by taking $p_1 \nearrow 1$.
If instead $d \in (-1,0)$ (i.e.\ Player Two
is better than Player One),
then $A(p_1,p_1-d)$ is maximised when
${d \over dp_1} \, A(p_1,p_1-d) = 0$,
i.e.\ $(1-p_1)^2+d = 0$,
i.e.\ $p_1 = 1 - \sqrt{-d} \in (0,1)$,
so there is a non-degenerate optimal choice of $p_1$.

Note, though, that
$1 - \sqrt{-d} > (1+d)/2$ for all $d \in (-1,0)$.
So, regardless of whether $d$ is positive or negative,
the Pops Rules still encourage more difficult shots than the
Traditional Rules do.
For example,
if $d=0.2$ so $p_2=p_1-0.2$,
then the optimal $p_1 \nearrow 1$
under TR, but the optimal $p_1 = (1+d)/2 = (1+0.2)/2 = 0.60$ under PR.
Or, if $d=-0.2$ so $p_2=p_1+0.2$,
then the optimal $p_1 = 1-\sqrt{-d}
= 1-\sqrt{0.2} \doteq 0.55$
under TR, but the optimal $p_1 = (1+d)/2 = (1-0.2)/2 = 0.40$ under PR.

\checkpage{2cm}
\bigskip\noindent\bf Ratio case: \rm

Next, consider the ratio case where $p_2=p_1/r$ for some fixed $r>0$.
To ensure that $p_1,p_2\in(0,1)$, we restrict to $p_1 \in (0, \min[1,r])$.
Then under PR,
\[
B(p_1,p_1/r)
\ = \ p_1\left(1 - \frac{p_1}{r}\right)
\ = \ p_1 - \frac{p_1^2}{r}
\, .
\]
If $r \geq 2$,
then $B(p_1,p_1/r)$ is an increasing function on its domain,
so it is optimal to choose $p_1 \nearrow 1$.
If instead $r<2$, then
$B(p_1,p_1/r)$ is maximised when $p_1 = r/2$, so $p_2=p_1/r=1/2$, i.e.\ the
optimal strategy for Player One (whether they are better or worse than
Player Two) is to take a shot that Player Two will make half the time.  

However, under TR,
\[
A(p_1,p_1/r)
\ = \ \frac{p_1 \left(1 - \frac{p_1}{r}\right)}{1 - p_1\left(\frac{p_1}{r}\right)}
\ = \ \frac{p_1 - \frac{p_1^2}{r}}{1 - \frac{p_1^2}{r}}
\, .
\]
To analyse this, we compute that for fixed $r$,
$$
{d \over dp_1} \, A(p_1,p_1/r)
\ = \ {r (p_1^2-2p_1+r) \over (p_1^2-r)^2}
\, .
$$
If $r > 1$ (i.e.\ Player One is better than Player Two), then
$r (p_1^2-2p_1+r) > r (p_1^2-2p_1+1) = r(p_1-1)^2 > 0$,
so ${d \over dp_1} \, A(p_1,p_1/r) > 0$ for all $p_1\in(0,1)$, i.e.\
$A(p_1,p_1/r)$ is an increasing function of $p_1$,
so it is again optimal to choose $p_1 \nearrow 1$.
If instead $r \in (0,1)$ (i.e.\ Player Two is better
than Player One), then $A(p_1,p_1/r)$ has critical points where
$p_1^2-2p_1+r=0$, i.e.\ from the quadratic formula
$p_1 = \Big[-(-2) \pm \sqrt{(-2)^2-4(1)(r)}\Big] \Big/ 2(1)
= 1 \pm \sqrt{1-r}$, and it follows that
$A(p_1,p_1/r)$ is maximised when $p_1 = 1 - \sqrt {1-r} \in (0,1)$.

Note, though, that $1 - \sqrt{1-r} > \frac{r}{2}$ for $r \in (0,1)$.
So, regardless of whether $r>1$ or $r<1$, the Pops Rules
encourage more difficult shots than the Traditional Rules do.
For example,
if $r=1.1$ so $p_2=p_1/1.1$,
then the optimal $p_1 \nearrow 1$
under TR, but the optimal $p_1 = 1.1/2 = 0.55$ under PR.
Or, if $r=0.8$ so $p_2=p_1/0.8$, then the optimal $p_1 = 1-\sqrt{1-r} =
1-\sqrt{0.2} \doteq 0.55$ under TR, but the optimal $p_1 = r/2 = 0.8/2 =
0.40$ under PR.

\bigskip

Our results can be summarised as follows:

\bigskip

\hskip -1cm
\begin{tabular}{|c|c|c|c|c|}
\hline
Case & $A(p_1,p_2)$ & opt.\ $p_1$ [TR]
& $B(p_1,p_2)$ & opt.\ $p_1$ [PR] \cr
\hline
\hline
$p_2=p_1$
  & ${p_1 \over 1+p_1}$ & $\nearrow 1$ & $p_1(1-p_1)$ & 1/2 \cr
\hline
$p_2=p_1-d \ \ (d \ge 0)$
  & ${p_1 - p_1^2 + p_1 d \over 1 - p_1^2 + p_1 d}$ & $\nearrow 1$
  & $(1+d) p_1 - p_1^2$ & $(1+d)/2$ \cr
\hline
$p_2=p_1-d \ \ (d < 0)$
  & ${p_1 - p_1^2 + p_1 d \over 1 - p_1^2 + p_1 d}$ & $1 - \sqrt{-d}$
  & $(1+d) p_1 - p_1^2$ & $(1+d)/2$ \cr
\hline
$p_2=p_1/r \ \ (r \ge 2)$
  & ${p_1 - {p_1^2 \over r} \over 1 - {p_1^2 \over r}}$ & $\nearrow 1$
  & $p_1 - {p_1^2 \over r}$ & $\nearrow 1$ \cr
\hline
$p_2=p_1/r \ \ (1 \le r < 2)$
  & ${p_1 - {p_1^2 \over r} \over 1 - {p_1^2 \over r}}$ & $\nearrow 1$
  & $p_1 - {p_1^2 \over r}$ & $r/2$ \cr
\hline
$p_2=p_1/r \ \ (0 < r < 1)$
  & ${p_1 - {p_1^2 \over r} \over 1 - {p_1^2 \over r}}$ & $1 - \sqrt{1-r}$
  & $p_1 - {p_1^2 \over r}$ & $r/2$ \cr
\hline
\end{tabular}

\bigskip\noindent
In particular, in both the linear and ratio cases, the following holds.  If
Player One is better than Player Two, then it is optimal under TR to
choose $p_1 \nearrow 1$.  But if Player Two is better than Player One, then
the optimal $p_1$ under TR is somewhere in $(0,1)$.  Nevertheless, in all
situations, the Pops Rules always lead to a smaller optimal $p_1$,
corresponding to more difficult (and, we believe, more interesting) shots
when compared to the Traditional Rules.

\bigskip

In light of this mathematical investigation into Horse probabilities and
optimisations, we feel that the Pops Rules provide a more interesting
alternative to the Traditional Rules, and should be used whenever playing
Horse on basketball courts throughout the land.

\bigskip\bigskip\noindent\bf Acknowledgements. \rm
We thank the reviewers and editor for helpful comments and
suggestions which have significantly improved the manuscript.

\vfil\eject
\section*{References}

NBA.com (2020),
All-Stars Chris Paul and Trae Young Headline First-Ever NBA HORSE
Challenge.  April 9, 2020.  Available at:
\hfil\break
https://www.NBA.com/news/nba-air-horse-challenge-espn

T.~Tarpey and R.T.~Ogden (2016), Statistical Modeling to Inform
Optimal Game Strategy: Markov Plays H-O-R-S-E".
The American Statistician {\bf 70(2)}, 181--186.

J.~Stevenson (2020), Shot Selection Tactics in the Basketball Game HORSE
and a Match between Mitchel Robinson and Brad Wanamaker.  Preprint.
Available at: http://dx.doi.org/10.2139/ssrn.3740706

\end{document}